\documentclass[a4paper,oneside,reqno,12pt]{amsart}
\usepackage{latexsym,amsthm,amscd,amsmath,amssymb}
\usepackage[mathscr]{eucal}

\advance\voffset 0cm \advance\hoffset 0cm

\def\n{\noindent}
\def\te{\text}

\def\vap{\varphi}
\def\we{\wedge}
\def\O{\Omega}

\newcommand{\R}{\mathbb R}
\newcommand{\C}{\mathbb C}

\newcommand{\B}{\mathbb B}

 \large 
\begin{document}
\title{On some weighted energy classes of plurisubharmonic functions}
\author{Le Mau Hai and Pham Hoang Hiep}
\maketitle
\begin{abstract} In  this paper we study the relation between the weighted energy class $\mathcal{E}_{\chi}$ introduced by S. Benelkouchi, V. Guedj and A. Zeriahi recently with the classes $\mathcal{E}$  and $\mathcal{N}$ studied by Cegrell. Moreover, we establish a generalized comparison principle for the operator $\text{M}_{\chi}$ and, as an application, we prove  a slight version of existence of solutions of Monge-Amp\`ere type equation in the class $\mathcal{E}_{\chi}(H,\O)$.
\vskip0.1cm

\vskip0.1cm
\n Keywords and Phrases:  class $\mathcal{E}_{\chi}$, class $\mathcal{E}$, class $\mathcal{N}$, strong comparison principle, operator $\text{M}_{\chi}$, Monge - Amp\`ere type equation.
\end{abstract}
\section{Introduction.}

\n To investigate  weighted energy classes of plurisubharmonic functions on a compact Kahl\"er manifold, as well as, in a hyperconvex domain $\O$ in $\C^n$ has been studied by several  authors recently. In the paper " The weighted Monge - Amp\`ere energy of quasiplurisubharmonic functions" V.Guedj and A.Zeriahi introduced and investigated the weighted energy class $\mathcal{E}_{\chi}(X,\omega)$ on a compact Kahl\"er manifold $X$ with a local potential $\omega$. They proved the continuity of weighted Monge - Amp\`ere operators on decreasing sequences of quasiplurisubharmonic functions in this class (see Theorem 2.6 in [GZ]) and investigated solutions of the complex Monge - Amp\`ere equation in this class (Theorem 4.1 there). In 2008,  Benelkouchi,  Guedj and Zeriahi introduced and studied the weighted energy class $\mathcal{E}_{\chi}$ as follows. Let $\O$ be a hyperconvex domain in $\C^n$ and  $\chi:(-\infty,0]\to [0,+\infty)$ be a decreasing function. Let
$$
\mathcal{E}_{\chi}=\mathcal{E}_{\chi}(\O)=\{u\in\text{PSH}^{-}(\O): \exists \   \ \mathcal{E}_0(\O)\ni u_j\searrow u,$$
 $$\sup\limits_{j}\int\limits_{\O}\chi(u_j)(dd^cu_j)^n < +\infty\}.
$$

\n (see [BGZ]). In Proposition 3.2 in their paper they proved that if $\chi(-\infty)= +\infty, \chi(0)\ne 0$ then $\mathcal{E}_{\chi}(\O)\subset\mathcal{F}^a(\O)$.
However, the relation between the classes $\mathcal{E}_{\chi}(\O)$ and $\mathcal{E}(\O)$ studied by Cegrell in [Ce3]  in [BGZ] is not clarified. Moreover, in [BGZ] there are no results about the relation between the classes $\mathcal{E}_{\chi}(\O)$ and $\mathcal{N}(\O)$ introduced in [Ce4]. Hence, the first aim of our paper is to find the relation between these classes. Namely, in section 3 in Theorem 3.1, relying on Theorem C of Kolodziej  in [Ko1],  we show that under certain conditions the class $\mathcal{E}_{\chi}(\O)$ is contained in $\mathcal{E}(\O)$ and $\mathcal{N}(\O)$. It should be remarked that our Theorem 3.1 is more general then Proposition A in [Bel] and the class $\mathcal{E}_{\chi}(\O)$ in the case $\chi(t)= \min(|t|^p, 1),  t<0$ has been introduced earlier in [H3].

\n The second aim  of the paper is to prove the strong principle for the operator $\text{M}_{\chi}$. Let $\chi:\R^{-}\times\O\to \R^{+}$ be such that $\chi(.,z)$ is a decreasing function. For each $u\in\mathcal{E}(\O)$ we put
$$\text{M}_{\chi}(u)=\chi(u(z),z)(dd^cu)^n.$$

\n  In the fourth section we prove the strong comparison principle for the operator $\text{M}_{\chi}$. This strong comparison principle will be used to study the operator $\text{M}_{\chi}$ instead of the complex Monge - Amp\`ere operator $(dd^cu)^n$. Consequently, we get the strong comparison principle of Xing stated and proved by Xing recently (see [Xi2]). Moreover, as an application of the strong comparison principle, we prove a slight version of existence of solutions of the equation of complex  Monge - Amp\`ere type in the class $\mathcal{E}_{\chi}(H,\O)$ which is similar to a recent result of Czyz in [Cz].
\vskip0.2cm
\n The paper is organized as follows. Besides the introduction, the paper has the three sections. In the second we recall some classes of  plurisubharmonic functions introduced and investigated by Cegrell and some authors recently. In particular, we give the definition of  the class $\mathcal{E}_{\chi}$ and the class $\mathcal{N}$. The third is devoted to establish the relation between the class $\mathcal{E}_{\chi}$ with the classes $\mathcal{E}$ and $\mathcal{N}$. The fourth says about the strong comparison principle for the operator $\text{M}_{\chi}$ and applications of this principle.

\section{Preliminaries}

\n{\bf 2.1.} We assume that readers are familiar with plurisubharmonic functions and the complex Monge - Amp\`ere operator for locally bounded  plurisubharmonic functions. Readers  can find notions about these objects in the monograph of Klimek [Kli] and the excellent paper of E.Bedford and B.A.Taylor [BT2].

\n{\bf 2.2.} Now we recall some classes of plurisubharmonic functions introduced and investigated by Cegrell recently (see [Ce2] and [Ce3]).

\n Let $\O\subset\C^n$ be a bounded hyperconvex domain, i.e there exists  a negative plurisubharmonic function $\varrho$ such that $\{z\in\O: \varrho(z)< -c\}\Subset\O$ for all $c>0$. By $\text{PSH}^{-}(\O)$ we denote the set of negative plurisubharmonic functions on $\O$. Following Cegrell we define the following.
\begin{align*}
\mathcal{E}_0=\mathcal{E}_0(\Omega)=\{\varphi\in\text{PSH}^{-}(\Omega)\cap\text{L}^{\infty}(\Omega):
\underset{z\to\partial{\Omega}}
\lim\varphi(z)=0,\
\ \int\limits_{\Omega}(dd^c\varphi)^n <\infty\}\\
\end{align*}
\n Now  for each $p>0$, put
\begin{align*}
\mathcal{E}_p=\mathcal{E}_p(\Omega)=\big\{\varphi\in\text{PSH}^{-}(\Omega) :\exists\   \  \mathcal{E}_0\ni\varphi_j\searrow\varphi,\   \  \underset{j}\sup\int\limits_{\Omega}(-\vap_j)^p(dd^c\varphi_j)^n<\infty\big\},\\
\end{align*}

\n
\begin{align*}
\mathcal{F}=\mathcal{F}(\Omega)=\big\{\varphi\in\text{PSH}^{-}(\Omega) :\exists\   \  \mathcal{E}_0\ni\varphi_j\searrow\varphi,\   \  \underset{j}\sup\int\limits_{\Omega}(dd^c\varphi_j)^n<\infty\big\}\\
\end{align*}
\n and
\begin{align*}
\mathcal{E}=\mathcal{E}(\Omega)=\big\{\varphi\in\text{PSH}^{-}(\Omega):\forall z_0\in\Omega, \exists   \  \ \text{a neighbourhood}\  \  \omega\ni z_0,\\
\mathcal{E}_0\ni\varphi_j\searrow\varphi \     \  \text{on}\  \ \omega,   \underset{j}\sup\int\limits_{\Omega}(dd^c\varphi_j)^n <\infty\big\}.
\end{align*}
The following inclusions  are  clear:  $\mathcal{E}_0\subset\mathcal{E}_p\subset\mathcal{E}$ and  $\mathcal{E}_0\subset\mathcal{F}\subset\mathcal{E}$.

\n In [Ce3] Cegrell has proved that the class $\mathcal{E}$ is the biggest on which the complex Monge - Amp\`ere operator $( dd^c.)^n$ exists. Moreover, if $u\in\mathcal{E}$ and $K\Subset\O$ then we can find $u_K\in\mathcal{F}$ such that $u= u_K$ on $K$. Another result of Cegrell in [Ce3] which we shall use  is following. Let $\{u_j, u\}\subset\mathcal{E}$ be such that $(dd^cu_j)^n$ is weakly convergent to $(dd^cu)^n$. If $g$ is a bounded plurisubharmonic function  then $g(dd^cu_j)^n$ is also weakly convergent to $g(dd^cu)^n$. We also use the following notation. Assume $\mathcal{K}\subset\mathcal{E}$. Then by $\mathcal{K}^a=\mathcal{K}^a(\O)$ we denote the set
$$\mathcal{K}^a=\mathcal{K}^a(\O)=\{ u\in\mathcal{K} : \   \  (dd^cu)^n << \text{Cap}\}.$$

\n That means that for every pluripolar set $E\subset\O$ we have $(dd^c u)^n (E) =0$.

\n{\bf 2.3.} We recall the class $\mathcal{N}(\O)$ introduced and investigated in [Ce4]. Let $\O$ be a hyperconvex domain in $\C^n$ and $\{\O_j\}_{j\geq 1}$ a fundamental sequence of $\O$. This is an increasing sequence of strictly pseudoconvex subsets $\O_j$ of $\O$ such that $\O_j\Subset \O_{j+1}$ and $\bigcup\limits_{j=1}^\infty\O_j = \O$. Let $\varphi\in\text{PSH}^{-}(\O)$. For each $j\geq 1$, put
$$\varphi^j= \sup\{ u: u\in\text{PSH}(\O),\     \  u\leq \varphi\    \  \text{on}\   \   \O\setminus \bar{\O}_j\}.$$

\n As in [Ce4], the function $\widetilde{\varphi}=\Bigl(\lim\limits_{j\to\infty}\varphi^j\Bigl)^*\in\text{PSH}(\O)$ and $\widetilde{\varphi}\in\mathcal{MPSH}(\O)$, where $\mathcal{MPSH}(\O)$ denotes the set of maximal plurisubharmonic functions on $\O$. Set
$$\mathcal{N}=\mathcal{N}(\O)=\{\varphi\in\mathcal{E}: \widetilde{\varphi}= 0\}$$

\n or equivalently,
$$\mathcal{N}=\mathcal{N}(\O)=\{\varphi\in\text{PSH}^{-}(\O): \varphi^j\uparrow 0\}.$$

\n It is easy to see that $\mathcal{E}_p\cup\mathcal{F}\subset\mathcal{N}$. Through this paper we assume that $H\in\mathcal{E}\cap \mathcal{MPSH}(\O)$. By $\mathcal{E}(H)=\mathcal{E}(H,\O)$(resp. $\mathcal{N}(H)=\mathcal{N}(H,\O)$) we denote the set
$$\mathcal{E}(H,\O)=\{ u\in\text{PSH}^{-}(\O): \exists \    \   \varphi\in\mathcal{E}\    \  \text{such that}\   \  H\geq u\geq \varphi + H\}$$
(resp. there exists $\varphi \in\mathcal{N}$ such that $H\geq u\geq \varphi + H$).

\n {\bf 2.4.} Now we recall  the class $\mathcal{E}_{\chi}(\O)$ introduced and investigated by Benelkourchi, Guedj and Zeriahi [BGZ] recently.

\n  Let $\chi: (-\infty,0]\to  [0,+\infty)$ be  a decreasing function. By $\mathcal{E}_{\chi}(\O)$ we denote the set of negative plurisubharmonic functions $u\in\text{PSH}^{-}(\O)$ such that there exists a sequence $\{u_j \}\subset\mathcal{E}_0(\O)$ with $u_j\searrow u$ and
 $$\sup\limits_{j}\int\limits_{\O}\chi(u_j)(dd^cu_j)^n <+\infty.$$

\n In the case $\chi(t)=  (-t)^p$(resp.$\chi(t)$ is bounded), $\mathcal{E}_{\chi}$ is the class $\mathcal{E}_p(\O)$(resp.$\mathcal{F}(\O)$) studied by Cegrell in [Ce2] and [Ce3].

\section {Relation between the classes $\mathcal{E}_{\chi}$, $\mathcal{E}$ and $\mathcal{N}$.}

\n The aim of this section is to establish the relation between the classes  $\mathcal{E}_{\chi}$, $\mathcal{E}$ and $\mathcal{N}$. Namely we prove the following.
\vskip0.1cm
\n{\bf Theorem 3.1.} {\sl Let $\chi:(-\infty,0]\to [0,+\infty)$ be a decreasing function satisfying $\chi\not\equiv 0$. Assume that $\{\varphi_j\}_{j\geq 1}\subset\mathcal{E}_0(\O)$ is such that
$$\sup\limits_{j}\int\limits_{\O}\chi(\varphi_j)(dd^c\varphi_j)^n <+\infty.$$

\n Then

a) $\varphi=\Bigl(\lim\sup\limits_{j }\varphi_j\Bigl)^*\in\mathcal{E}(\O)$.

b) Moreover, if $\chi(t)> 0$ when $t<0$ then $\varphi=\Bigl(\lim\sup\limits_{j }\varphi_j\Bigl)^*\in\mathcal{N}(\O)$}

\n In order to prove the above theorem we need the following.
\vskip0.1cm
\n {\bf Lemma 3.2.}{\sl  Let $\chi:(-\infty,0]\to[0,+\infty)$ be a decreasing function and $u,v\in\mathcal{E}_0(\O)$ satisfy $\chi(u)(dd^cu)^n = (dd^cv)^n$. Then
$$ u\geq \frac{1}{\sqrt[n]{\chi(-t)}}v - t,\     \  \forall \   \ t>0 \    \ \text{such that}\    \ \chi(-t) > 0.$$}
\vskip0.1cm
\n{\it Proof.} It is easy to see that $u\geq -t\geq \frac{1}{\sqrt[n]{\chi(-t)}}v - t$ on the set $\{u\geq -t\}$. It suffices to show  that $u\geq w=\frac{1}{\sqrt[n]{\chi(-t)}}v - t$ on the open set $D=\{u<-t\}$. We have that
$$ \lim\limits_{z\to\partial{D}}u(z)=-t\geq \lim\limits_{z\to\partial{D}}w(z).$$

\n We show that $(dd^cu)^n\leq (dd^cw)^n$ on $D=\{u<-t\}$. Indeed, since $\chi$ is decreasing then $\chi(-t)(dd^cu)^n\leq\chi(u)(dd^cu)^n=(dd^cv)^n$ on $D=\{u<-t\}$. Hence, $(dd^cu)^n\leq\frac{1}{\chi(-t)}(dd^cv)^n= (dd^cw)^n$ on $D=\{u<-t\}$. Using the comparision principle for $u,w$ on the set $D=\{u<-t\}$ it follows that $u\geq w $ on $D$ and the desired conclusion of the lemma follows.
\vskip0.1cm
\n {\bf Proof of Theorem 3.1.} By theorem C in [Ko1] we can find $v_j\in\mathcal{E}_0(\O)$ such that $\chi(\varphi_j)(dd^c\varphi_j)^n = (dd^cv_j)^n$. Lemma 3.2 implies that $\varphi_j\geq \frac{1}{\sqrt[n]{\chi(-t)}}v_j - t$ for all $t>0$ such that $\chi(-t)>0$. From
$$\sup\limits_{j }\int\limits_{\O}(dd^cv_j)^n= \sup\limits_{j }\int\limits_{\O}\chi(\varphi_j)(dd^c\varphi_j)^n <+\infty$$

\n then Theorem 3.7 in [KH] implies that $v=\Bigl(\lim\sup\limits_{j } v_j\Bigl)^*\in\mathcal{F}(\O)$. Hence, $\varphi=\Bigl(\lim\sup\limits_{j }\varphi_j\Bigl)^*\geq \frac{1}{\sqrt[n]{\chi(-t)}}v-t$ for all $t>0$ such that $\chi(-t)>0$. Since $\chi\not\equiv 0$ then  we can  choose $t>0$ such that $\chi(-t)> 0$. Then $\varphi\geq \frac{1}{\sqrt[n]{\chi(-t)}}v-t$ and from  $v\in\mathcal{F}(\O)$ it follows that $\varphi\in\mathcal{E}(\O)$ and a) in Theorem 3.1 is proved.

\n  Now we give the proof of b) of Theorem 3.1. From $v\in\mathcal{F}(\O)$ we infer that $\widetilde{v}=0$. Thus $\widetilde{\varphi}\geq -t$. Assume that $\chi(-t)>0$ for all $t>0$. Then  $\widetilde{\varphi}\geq -t$ for all $t>0$. Hence, $\widetilde{\varphi}=0$ and it follows that $\varphi\in\mathcal{N}(\O)$. Theorem 3.1 is completely proved.
\vskip0.1cm
\n From Theorem 3.1 we derive the following corollary.
\vskip0.1cm
\n {\bf Corollary 3.3.} a) {\sl $\mathcal{E}_{\chi}(\O)\subset\mathcal{E}(\O)$ if $\chi$ is decreasing and $\chi\not\equiv 0$.}

b) {\sl $\mathcal{E}_{\chi}(\O)\subset\mathcal{N}(\O)$ if $\chi$ is decreasing and $\chi(t)>0$ for all $t<0$.}

\n We have some the following remarks concerning to Theorem 3.1.

\n{\bf Remark 3.4.}

{\bf (i)} $\chi(-\infty)=+\infty$ if and only if  $\mathcal{E}_{\chi}\subset\mathcal{E}^a$.

\n Indeed, first assume that $\chi(-\infty)=+\infty$ and $\varphi\in\mathcal{E}_{\chi}$. There exists a sequence $\{\varphi_j\}\subset\mathcal{E}_0$ such that $\varphi_j\searrow \varphi$ and
$$ \sup\limits_{j}\int\limits_{\O}\chi(\varphi_j)(dd^c\varphi_j)^n < +\infty.$$

\n Since $\chi$ is decreasing then we have
\begin{align*}
&\int\limits_{\{\varphi_j<-t\}}(dd^c\varphi_j)^n\leq\int\limits_{\{\varphi_j<-t\}}\frac{\chi(\varphi_j)}{\chi(-t)}(dd^c\varphi_j)^n\\
& \leq \frac{\sup\limits_{j}\int\limits_{\O}\chi(\varphi_j)(dd^c\varphi_j)^n}{\chi(-t)}.
\end{align*}

\n Since $\{\varphi_j<-t\}$ is an increasing sequence of subsets which tends to $\{\varphi< -t\}$ then by tending $j\to\infty$ we get
$$\int\limits_{\{\varphi<-t\}}(dd^c\varphi)^n\leq \frac{\sup\limits_{j}\int\limits_{\O}\chi(\varphi_j)(dd^c\varphi_j)^n}{\chi(-t)}.$$

\n By letting $t\to +\infty$ we infer that
$$\int\limits_{\{\varphi=-\infty\}}(dd^c\varphi)^n=0.$$

\n Hence, $\varphi\in\mathcal{E}^a$.

\n Second, let $\chi(-\infty)<+\infty$. Then  it follows that $\mathcal{F}\subset\mathcal{E}_{\chi}$. However, $\mathcal{F}\not\subset\mathcal{E}^a$ then $\mathcal{E}_{\chi}$ is not a subset of $\mathcal{E}^a$.

{\bf(ii)} The condition {\bf b)} in our Theorem 3.1 is sharp. Indeed, assume that there exists $t_0> 0$ such that $\chi(-t_0)=0$. Then the function $\varphi = -t_0\in\mathcal{E}_{\chi}$ because we have
$$\int\limits_{\O}\chi(\varphi_j)(dd^c\varphi_j)^n = 0$$

\n for all $0\geq \varphi_j > -t_0$. But $\varphi=-t_0\notin\mathcal{N}$.

\section{ Strong comparison principle for the operator $\text{M}_{\chi}$.}

\n As be said above, the aim of this section is to give the strong comparison principle for the operator $\text{M}_{\chi}$ and to find applications of this principle. Through the section we assume that $\chi:(-\infty,0]\to[0,+\infty)$ is a decreasing function. For each a such  function $\chi$ we set
$$\chi_k (t) = \int\limits_t^0 dt_1\int\limits_{t_1}^0 dt_2 ... \int\limits_{t_{k-1}}^0 \chi (t_k)dt_k.$$
 First we prove the following.
\vskip0.2cm
\n
{\bf Theorem 4.1.}  {\it Let $u\in\mathcal E^a(\O)$ and $v\in\mathcal{E}(\O)$ be such that $\liminf\limits_{z\to\partial\Omega}[u(z)-v(z)]\geq 0$. Then we have
$$\int\limits_{\{u<v\}} \chi_k(u-v) dd^c w_1\we ...\we dd^c w_k\we T$$
$$\leq \int\limits_{\{u<v\}} -w_1 \chi (u-v) \Bigl[(dd^cu)^k-(dd^cv)^k\Bigl]\we T$$
for every $w_1,...,w_k\in\te{PSH}^-(\O)\cap\te{L}^\infty(\O)$, $w_2,...,w_k\geq -1$ and $T=dd^c w_{k+1}\we ...\we dd^cw_n$ with $w_{k+1},...,w_n\in\mathcal E(\Omega)$.}
\vskip0.1cm
\n In order to prove the above theorem we need the following.
\vskip0.2cm
\n
{\bf Lemma 4.2.}  {\it Let $u,v\in\te{PSH}(\O)\cap\te{L}_{loc}^\infty(\O)$ be such that $u=v$ on $\O\setminus K$ with $K\Subset\O$ and $u\leq v$. Then we have
$$\int\limits_{\O} \chi_1(u-v) dd^cw\we T\leq \int\limits_{\O} -w \chi (u-v) dd^c(u-v)\we T,$$
for every $w\in\te{PSH}^-(\O)\cap\te{L}^\infty(\O)$ and $T= dd^c w_1\wedge\cdots\wedge dd^cw_{n-1}$, where $ w_1,\cdots, w_{n-1}\in\mathcal{E}$. }
\vskip0.2cm
\n
{\it Proof.}  First, we consider the case $\chi\in C^1(\mathbb R)$, $u,v\in C^\infty (\Omega )$. By Stoke's theorem we have
\begin{align*}
&\int\limits_{\O}\chi_1 (u-v) dd^c w \we T=\int\limits_{\O} w dd^c \chi_1 (u-v)\we T\\
&=\int\limits_{\O}w d\Bigl(-\chi (u-v) d^c(u-v)\Bigl)\we T\\
&=\int\limits_{\O}-w \chi '(u-v)d(u-v)\we d^c(u-v)\we T + \int\limits_{\O} -w \chi (u-v) dd^c(u-v)\we T\\
&\leq \int\limits_{\O}-w \chi (u-v) dd^c(u-v)\we T.
\end{align*}

\n
Second, we consider the case $\chi\in C^1(\mathbb R)$. Set $u_\delta =u*\rho_\delta$, $v_\delta =v*\rho_\delta$ on $\Omega_\delta = \{z\in\Omega: d(z,\mathbb C^n\backslash\Omega)>\delta\}$, here $\rho_\delta (z) =\frac 1 {\delta^{2n}} \rho (\frac z {\delta})$ with $\rho\in C_0^\infty (\mathbb C^n)$ such that supp$\rho\subset \B=\{z\in\mathbb C^n:\ ||z||<1\}$, $\rho (z)=\rho (|z|)\geq 0$ and $\int\limits_{\B}\rho (z) dV_{2n} (z) =1$. We can find $\delta_0>0$ such that $u_\delta =v_\delta$ on $\Omega_\delta\setminus K_{\delta}$ for all $\delta<\delta_0$, where $K_{\delta}=\{z\in\O: d(z,K)\leq \delta\}$. By the  first case we get
$$\int\limits_{\O_\delta} \chi_1(u_\delta-v_\delta) dd^cw\we T\leq \int\limits_{\O_\delta} -w \chi (u_\delta-v_\delta) dd^c(u_\delta-v_\delta)\we T.$$
 The quasicontinuity of $u,v$ (Theorem 3.5 in [BT2]) together with  tending  $\delta\to 0$ implies that
$$\int\limits_{\O} \chi_1(u-v) dd^cw\we T\leq \int\limits_{\O} -w \chi (u-v) dd^c(u-v)\we T.$$

\n
Third, we consider the case $\chi (-\infty)<+\infty$. We find decreasing functions $\chi_j\in C^1(\mathbb R)$ such that $\chi_j\to\chi$ on $(-\infty, 0]$. By the second case we have
$$\int\limits_{\O} \chi_{j1}(u-v) dd^cw\we T\leq \int\limits_{\O} -w \chi_j (u-v) dd^c(u-v)\we T.$$
Using Lebesgue monotone  convergence theorem, by letting $j\to\infty$ we get
$$\int\limits_{\O} \chi_1(u-v) dd^cw\we T\leq \int\limits_{\O} -w \chi (u-v) dd^c(u-v)\we T.$$
In the  general case, set $\chi_j=\min (\chi, j)$. By the third case we have
$$\int\limits_{\O} \chi_{j1}(u-v) dd^cw\we T\leq \int\limits_{\O} -w \chi_j (u-v) dd^c(u-v)\we T.$$
 Lebesgue  monotone convergence theorem implies that
$$\int\limits_{\O} \chi_1(u-v) dd^cw\we T\leq \int\limits_{\O} -w \chi (u-v) dd^c(u-v)\we T.$$
The proof of the lemma is complete.
\vskip0.1cm
\n
{\bf Lemma 4.3.}  {\it Let $u,v\in\mathcal E(\O)$ be such that $u=v$ on $\O\setminus K$ with $K\Subset\O$ and $u\leq v$. Then we have
\begin{align*} \int\limits_{\O} \chi_1(u-v) dd^cw\we T&\leq \int\limits_{\{v>-\infty\}} -w \chi (u-v) dd^c(u-v)\we T\\
&+\chi (-\infty)\int\limits_{\{u=-\infty\}} -w dd^cu\we T\\
&-\chi (0)\int\limits_{\{v=-\infty\}} -w dd^cv\we T,\hskip3.3cm (1)
\end{align*}

\n for every $w\in\te{PSH}^-(\O)\cap\te{L}^\infty(\O)$ and $T=dd^c w_1\wedge ...\wedge dd^c w_{n-1}$ with $w_1,...,w_{n-1}\in\mathcal E(\Omega)$.}
\vskip0.2cm
\n
{\sl Proof.} Choose $\phi\in C_0^\infty (\Omega)$ such that $0\leq\phi\leq 1$ and $\phi =1$ on a neighbourhood of $K$.

\n
First, we consider the case $\chi(-\infty )<+\infty$. Set
$$u_j=\max (u,-j),\ v_j=\max (v,-j).$$
By Lemma 4.2 we have
\begin{align*}& \int\limits_{\O} \chi_1(u_j-v_j) dd^cw\we T\leq \int\limits_{\O} -w \chi (u_j-v_j) dd^c(u_j-v_j)\we T\\
&=\int\limits_{\O} -w \phi\chi (u_j-v_j) dd^c(u_j-v_j)\we T.\hskip5.3cm(2)
\end{align*}
Since  $u\leq v$ then by simple  arguments we notice that  $u_j-v_j\searrow u-v$. Using Lebesgue  monotone convergence theorem we have
$$  \lim\limits_{j\to\infty}\int\limits_{\O} \chi_1(u_j-v_j) dd^cw\we T
=\int\limits_{\O} \chi_1(u-v) dd^cw\we T\hskip3.5cm(3)$$
By Theorem 4.1 in [KH], we have
\begin{align*}\lim\limits_{j\to\infty}\int\limits_{\O} -w\phi \chi (u_j-v_j) dd^cv_j\we T&=\lim\limits_{j\to\infty}\int\limits_{\{v>-j\}} -w\phi \chi (u_j-v_j) dd^cv\we T\\
&+\chi(0)\lim\limits_{j\to\infty}\int\limits_{\{v=-j\}} -w\phi dd^cv_j\we T\\
&=\int\limits_{\{v>-\infty\}} -w\phi \chi (u-v) dd^cv\we T\\
&+\chi(0)\lim\limits_{j\to\infty}\int\limits_{\Omega} -w\phi dd^cv_j\we T\\
&-\chi(0)\lim\limits_{j\to\infty}\int\limits_{\{v>-j\}} -w\phi dd^cv\we T\\
&=\int\limits_{\{v>-\infty\}} -w\phi \chi (u-v) dd^cv\we T\\
&+\chi(0)\int\limits_{\{v=-\infty\}}-w\phi dd^cv\we T.\hskip1cm(4)
\end{align*}
\n Now Theorem 4.1 in [KH] implies that
\begin{align*}& \limsup\limits_{j\to\infty}\int\limits_{\O} -w\phi \chi (u_j-v_j) dd^cu_j\we T\leq\lim\limits_{j\to\infty}\int\limits_{\{u>-j\}} -w\phi \chi (u_j-v_j) dd^c u\we T\\
&+\chi(-\infty )\lim\limits_{j\to\infty}\int\limits_{\{u=-j\}} -w\phi dd^cu_j\we T\\
&=\int\limits_{\{u>-\infty\}} -w\phi \chi (u-v) dd^cu\we T\\
&+\chi(-\infty )\int\limits_{\{u=-\infty\}}-w\phi dd^c u\we T.\hskip6.1cm(5)
\end{align*}
Put  (2), (3), (4) and (5) together and tend $\phi\nearrow 1$ we get  (1).

\n
In the general case we set $\chi_j = \min (\chi,j)$. Applying the  first case to $\chi_j$ and using Lebesgue  monotone convergence theorem we obtain (1).

\n
The proof of the lemma is complete.

\n
{\sl Proof of Theorem 4.1.}  Set $v_\epsilon =\max (u,v-\epsilon)\in\mathcal E^a (\Omega)$. Using Lemma 4.3 for $u,v_\epsilon$ we get
\begin{align*}
&\int\limits_{\O} \chi_k(u-v_\epsilon ) dd^cw_1\we ...\we dd^cw_k\we T\\
&\leq \int\limits_{\Omega} \chi_{k-1} (u-v_\epsilon ) dd^cu\we dd^cw_1\we ...\we dd^cw_{k-1}\we T\\
&\leq\cdots\cdots\cdots\cdots\cdots\cdots\cdots\cdots\cdots\\
&\leq \int\limits_{\Omega} \chi_1 (u-v_\epsilon ) dd^cw_1\we (dd^c u)^{k-1}\we T\\
&\leq \int\limits_{\Omega} \chi_1 (u-v_\epsilon ) dd^cw_1\we \Bigl[\sum\limits_{j=1}^{k}(dd^c u)^{k-j}\we (dd^c v_\epsilon)^{j-1}\Bigl]\we T\\
&\leq \int\limits_{\Omega } -w_1 \chi (u-v_\epsilon) \Bigl[(dd^c u)^k-(dd^c v_\epsilon)^k\Bigl]\we T\\
&\leq \int\limits_{\{u\leq v-\epsilon\}} -w_1 \chi (u-v+\epsilon) (dd^c u)^k\we T\\
&-\int\limits_{\{u<v-\epsilon\}} -w_1 \chi (u-v+\epsilon) (dd^c v)^k\we T.
\end{align*}
By letting  $\epsilon\to 0$ we obtain
$$\int\limits_{\{u<v\}} \chi_k(u-v) dd^c w_1\we ...\we dd^c w_k\we T\leq \int\limits_{\{u<v\}} -w_1 \chi (u-v) \Bigl[(dd^cu)^k-(dd^cv)^k\Bigl]\we T$$

\n and Theorem 4.1 is completely proved.

\n We have the following result for the case $u\in\mathcal{N}^a(H)$ and $v\in\mathcal{E}(H)$.
\vskip0.2cm
\n
{\bf Theorem 4.4.} {\it Let $u\in\mathcal N^a(H)$, $v\in\mathcal E(H)$. Then we have
$$\int\limits_{\{u<v\}} \chi_k(u-v) dd^c w_1\we ...\we dd^c w_k\we T\leq \int\limits_{\{u<v\}} -w_1 \chi (u-v) \Bigl[(dd^cu)^k-(dd^cv)^k\Bigl]\we T$$

\n for every $w_1,...,w_k\in\te{PSH}^-(\O)\cap\te{L}^\infty(\O)$, $w_2,...,w_k\geq -1$ and $T=dd^c w_{k+1}\we ...\we dd^cw_n$ with $w_{k+1},...,w_n\in\mathcal E(\Omega)$.}
\vskip0.2cm
\n
{\it Proof.} Take $\varphi\in\mathcal N^a$ such that $u\geq H+\varphi$. Since $(dd^c H)^n=0$ it follows that $u\in\mathcal{E}^a$. Set $v_j=\max (u, v+\varphi^j)$. Then $v_j\in\mathcal{E}$ and $\lim\inf\limits_{z\to\partial{\O}} (u(z)- v_j(z))\geq 0$. Applying Theorem 4.1 to $u, v_j$ we have
$$\int\limits_{\{u<v_j\}} \chi_k(u-v_j) dd^c w_1\we ...\we dd^c w_k\we T\leq \int\limits_{\{u<v_j\}} -w_1 \chi (u-v_j) \Bigl[(dd^cu)^k-(dd^cv_j)^k\Bigl]\we T$$

\n From $\{u<v_j\}=\{u<v+\varphi^j\}$ and by Theorem 4.1 in [KH] we get
$$
\int\limits_{\{u<v+\varphi^j\}} \chi_k(u-v-\varphi^j) dd^c w_1\we ...\we dd^c w_k\we T$$
$$\leq \int\limits_{\{u<v\}} -w_1 \chi (u-v-\varphi^j) \Bigl[(dd^cu)^k-(dd^cv)^k\Bigl]\we T$$

\n Tending  $j\to\infty$ we arrive that
$$\int\limits_{\{u<v\}} \chi_k(u-v) dd^c w_1\we ...\we dd^c w_k\we T\leq \int\limits_{\{u<v\}} -w_1 \chi (u-v) \Bigl[(dd^cu)^k-(dd^cv)^k\Bigl]\we T$$

\n and the desired conclusion follows.

\n Now we give some corollaries.
\vskip0.2cm
\n
{\bf Corollary 4.5.} {\it Let $u\in\mathcal N^a(H)$, $v\in\mathcal E(H)$. Then we have
$$\int\limits_{\{u<v\}} \chi_n(u-v) dd^c w_1\we ...\we dd^c w_n\leq \int\limits_{\{u<v\}} -w_1 \chi (u-v) \Bigl[(dd^cu)^n-(dd^cv)^n\Bigl]$$

\n for every $w_1,...,w_n\in\te{PSH}^-(\O)\cap\te{L}^\infty(\O)$, $w_2,...,w_n\geq -1$.}

\n
Now assume that $\chi (t) = |t|^p= (-t)^p, t<0$. Then we have $\chi_k(t)=\frac{1}{(p+k)\cdots ( p+1)}(-t)^{p+k}$ and we obtain the strong comparison principle  of Xing in [Xi2].
\vskip0.1cm
\n
{\bf Corollary 4.6.} {\it Let $u\in\mathcal E_p(\O)$, $v\in\mathcal E(\Omega )$. Then
$$\frac{1}{(n+p)\cdots(1+p)}\int\limits_{\{u<v\}}(v-u)^{n+p} dd^cw_1\we\cdots\we dd^cw_n$$
$$\leq\int\limits_{\{u<v\}} -w_1(v-u)^p\Bigl[ (dd^cu)^n - (dd^c v)^n\Bigl]$$

\n for all $w_1,\cdots, w_n\in\te{PSH}^-(\O)\cap\te{L}^\infty(\O), -1\leq w_2, \cdots, w_n\leq 0$.}
\vskip0.2cm
\n
{\bf Theorem 4.7.} {\sl Let $u\in\mathcal N^a(H)$, $v\in\mathcal E(H)$ be such that $(dd^c u)^n\leq (dd^c v)^n$ on $\{u<v\}$. Then $u\geq v$.}

\n
{\sl Proof.} Take $\varphi\in\mathcal N^a$ such that $u\geq H+\varphi$. We have $\{u<v+\varphi^j\}\Subset\Omega$ and $(dd^c u)^n\leq (dd^c( v+\varphi^j))^n$ on $\{u<v+\varphi^j\}$. By comparison principle we get $u\geq v+\varphi^j$. Letting $j\to\infty$ we obtain $u\geq v$.
\vskip0.2cm
\n In the end of the paper we prove a slight version of existence of solutions of the equation of complex Monge-Amp\`ere type equation in the class $\mathcal{E}_{\chi}(H)$. It should be noted that equations of complex Monge- Amp\`ere type were studied earlier by S.Kolodziej and U.Cegrell in [CeKo]. Recently R.Czyz has solved
this equation in the class $\mathcal{E}_{\chi}$ (see [Cz]). Here we give a simple proof of a slight version of existence of solutions of this equation which is similar to  [Cz].

\n For each $\chi:\R^{-}\times\O\longrightarrow \R^{+}$ we set
$$\te{M}_{\chi}(u) = \chi(u(z),z)(dd^cu)^n$$

\n where $u\in\mathcal{E}(\O)$.

\n First we prove the following comparison principle for the operator $\text{M}_{\chi}$.
\vskip0.1cm
\n{\bf Theorem 4.8.} {\it Let $\chi:\R^-\times\O\longrightarrow \R^{+}$ be such that $\chi(.,z)$ is decreasing for all $z\in\O$. Given $u\in\mathcal N^a(H)$, $v\in\mathcal E(H)$ such that $\te{M}_{\chi}(u)\leq \te{M}_{\chi}(v)$.

\n
Then $u\geq v$.}

\n
{\it Proof.} We have
$$(dd^c u)^n\leq \frac {\chi(v(z),z)}{\chi(u(z),z)} (dd^c v)^n\leq (dd^c v)^n,$$
on $\{u<v\}$. Theorem 4.7  implies that $u\geq v$.
\vskip0.2cm
\n
{\bf Corollary 4.9. }{\it Let $\chi_1,\chi_2: \R^-\times\O\longrightarrow \R^{+}$ be such that $\chi_1(.,z),\chi_2(.,z)$ are decreasing functions and $\chi_2\leq \chi_1$ on $\R^-\times\O$. Given $u_1\in\mathcal{N}^a(H)$, $u_2\in\mathcal E(H)$ such that $\te{M}_{\chi_1}(u_1)\leq \te{M}_{\chi_2}(u_2)$.

\n
Then $u_1\geq u_2$.}
\vskip0.2cm
\n{\it Proof.}  We have  $\te{M}_{\chi_1}(u_1)\leq\te{M}_{\chi_2}(u_2)\leq\te{M}_{\chi_1}(u_2)$. From Theorem 4.8 it follows that $u_1\geq u_2$.

\n Now let $\chi:\R^-\times\O\longrightarrow\R^{+}$ be such that $\chi(.,z)$ is a decreasing function for all $z\in\O$ and $\inf\limits_{z\in\O}\chi(t,z) > 0, \forall\  \  t<0$. Let $H\in\mathcal{E}(\O)\cap\mathcal{MPSH}(\O)$ be given. We define
\begin{align*}\mathcal{E}_{\chi}(H,\O)=\{\vap\in\mathcal{N}(H):\   \ \exists \  \ \mathcal{E}_0(H)\ni\vap_j\searrow\vap,\\
\sup\limits_{j\geq 1}\int\limits_{\O}\te{M}_{\chi}(\vap_j)<+\infty\}.
\end{align*}

\n Let $\mu$ be a non-negative finite measure on $\O$ and $\mu(P)=0$ for all pluripolar sets $P\subset\O$. We give a simple proof for a slight version of Main Theorem in [Cz].
\vskip0.1cm
\n{\bf Theorem 4.10.} {\it Let $\chi$ and $\mu$ satisfy all the above assumptions. Then there exists an unique function $u\in\mathcal{E}_{\chi}(H,\O)$ such that $\te{M}_{\chi}(u) = \mu$.}
\vskip0.1cm
\n{\it Proof.}  By [Ce3] we can find $\Phi\in\mathcal{E}_0(\O)$ and $0\leq f\in\te{L}^1_{loc}((dd^c\Phi)^n)$ such that
$$\mu = f(dd^c\Phi)^n.$$

\n Set $\mu_j=1_{\O_j}\min(f,j)(dd^c\Phi)^n$ where $\O_j\nearrow \O$ is defined as in section 2. Let $\eta_j\nearrow \frac{1}{\chi}=\eta$ such that $\eta_j(.,z)\nearrow$ for all $z\in\O$ and $\eta_j\in\te{C}^\infty(\R^-\times\O)$. Put $\chi_j=\frac{1}{\eta_j}\searrow \chi$. By [CeKo] we can find $u_j\in\mathcal{F}^a(H)$ such that
$$ (dd^cu_j)^n= \eta_jd\mu_j= \frac{d\mu_j}{\chi_j}.$$

\n Hence, $\te{M}_{\chi_j}(u_j)= d\mu_j$. Corollary 4.9 implies that $u_j\searrow u$. Now we show that $u_j\in\mathcal{E}_0(H)$. By [CeKo] we can find $\varphi_j\in\mathcal{F}^a$ such that
$$(dd^c\varphi_j)^n=\eta_j d\mu_j=\frac{d\mu_j}{\chi_j}.$$

\n Hence, $\text{M}_{\chi_j}(\varphi_j)= d\mu_j$. Note that $\eta_j(\varphi_j(z), z),  z\in\O_j$, is above bounded and
$$(dd^c\varphi_j)^n=\eta_j(\varphi_j(z),z) d\mu_j= \eta_j(\varphi_j(z),z) 1_{\O_j}\min(f,j)(dd^c\Phi)^n$$

\n then by the comparison principle it follows that $\varphi_j\in\mathcal{E}_0$. We prove that $u_j\geq H+\varphi_j$. Indeed, because $\{u_j, H+\varphi_j\}\subset\mathcal{F}^a(H)$ and
$$\text{M}_{\chi_j}(u_j)= d\mu_j= \text{M}_{\chi_j}(\varphi_j)\leq \text{M}_{\chi_j}(H+\varphi_j)$$

\n then Theorem 4.8 implies that $u_j\geq H+\varphi_j$. Moreover, by the above argument we infer that $\{\varphi_j\}$ is decreasing and we have
\begin{align*}
\sup\limits_{j\geq 1}\int\limits_{\O}\te{M}_{\chi}(\varphi_j)\leq \sup\limits_{j\geq 1}\int\limits_{\O}\te{M}_{\chi_j}(\varphi_j)\\
=\sup\limits_{j\geq 1}\int\limits_{\O}d\mu_j= \mu(\O).
\end{align*}

\n Now by $\inf\limits_{z\in\O}\chi(t,z)>0$ for all $t<0$ and using a similar argument as in the proof of  Theorem 3.1 we derive  that $\varphi=\lim\limits_{j\to\infty}\varphi_j\in\mathcal{N}$. Since $H\geq u_j\geq H+\varphi_j$ for all $j$ it follows that $H\geq u\geq H+\varphi$. Thus $u\in\mathcal{N}(H)$. At the same time, we have
\begin{align*}
\sup\limits_{j\geq 1}\int\limits_{\O}\te{M}_{\chi}(u_j)\leq \sup\limits_{j\geq 1}\int\limits_{\O}\te{M}_{\chi_j}(u_j)\\
=\sup\limits_{j\geq 1}\int\limits_{\O}d\mu_j= \mu(\O)
\end{align*}
\n and this shows that $u\in\mathcal{E}_{\chi}(H,\O)$. Now we have
$$\te{M}_{\chi}(u)=\lim\limits_{j\to\infty}\te{M}_{\chi_j}(u_j)=\lim\limits_{j\to\infty}d\mu_j = d\mu.$$

\n From Theorem 4.8 it follows the uniqueness of $u$ and the theorem is proved.

\n{\bf Acknowledgement.} Authors would like to thank Per Ahag and Slawomir Dinew  for valuable comments and refering us to the paper of  S. Benelkouchi.

\vskip0.2cm

\n Department of Mathematics

\noindent Hanoi National University of Education

\noindent Tuliem - Hanoi - Vietnam

\n Email : mauhai@ fpt.vn

\hskip1cm phhiep$_-$vn@yahoo.com
\end{document}